\documentstyle[12pt]{article} 
\newfont{\twlvmsb}{msbm10 scaled\magstep1}
\newfont{\ninemsb}{msbm9}
\newfont{\sixmsb}{msbm6}
\newfam\msbfam
\textfont\msbfam=\twlvmsb
\scriptfont\msbfam=\ninemsb
\scriptscriptfont\msbfam=\sixmsb
\catcode`\@=11
\def\Bbb{\ifmmode\let\next\Bbb@\else
  \def\next{\errmessage{Use \string\Bbb\space only in math mode}}\fi\next}
\def\Bbb@#1{{\Bbb@@{#1}}}
\def\Bbb@@#1{\fam\msbfam#1}
\newfont{\largeeufm}{eufm10 scaled\magstep4}
\newfont{\twlveufm}{eufm10 scaled\magstep1}
\newfont{\elveufm}{eufm10 at 11pt}
\newfont{\teneufm}{eufm10}
\newfont{\nineeufm}{eufm9}
\newfam\eufam
\textfont\eufam=\twlveufm
\scriptfont\eufam=\nineeufm
\def\frak{\ifmmode\let\next\frak@\else
\def\next{\errmessage{Use \string\frak\space only in math mode}}\fi\next}
\def\frak@#1{{\fam\eufam{{#1}}}}
\catcode`@=12
\newcommand{\g}{\widehat{sl}(2)} 
\newcommand{\U}{U(\widehat{sl}(2))} 
\newcommand{\Z}{{\Bbb Z}} 
\newcommand{\C}{{\Bbb C}} 
 
\newcommand{\ba}{\begin{eqnarray}}
\newcommand{\na}{\end{eqnarray}}
\newcommand{\ban}{\begin{eqnarray*}}
\newcommand{\nan}{\end{eqnarray*}}

\newcommand{\BPW}{Poincar\'e - Birkhoff - Witt\ } 
\newcommand{\BBW}{Bott - Borel - Weil\ } 

\begin{document} 
\title{\large Vector Coherent State Realization of 
Representations of \\
\large the Affine Lie Algebra $\widehat{sl}(2)$ } 
\author{\small R. B. ZHANG \\
\small Department of Pure Mathematics,  University of Adelaide, 
Adelaide, Australia} 
\date{} 
\maketitle

\begin{abstract}  
The method of vector coherent states is generalized 
to study representations of the affine Lie algebra $\widehat{sl}(2)$.  
A large class of highest weight irreps  
is explicitly constructed, which contains the integrable 
highest weight irreps as special cases. 
\end{abstract} 

\vspace{1cm}
\section{\normalsize INTRODUCTION} 
The method of vector coherent states was independently developed   
by the research groups of Quesne and Rowe\cite{Hecht}\cite{Rowe} 
to study representations of Lie groups appearing in physics. 
Applied to the irreducible representations of the 
compact semi - simple Lie groups, the method 
enables one to obtain explicit realizations of the generators 
of the corresponding 
Lie algebras in terms of holomorphic differential operators, 
and more importantly, the K - matrix technique \cite{Rowe}  
of the method allows the identification of 
the subsets of holomorphic polynomials 
which form the irreducible representation spaces.  
The method also provides a powerful machinery for constructing 
unitary representations of noncompact Lie groups; we refer to 
\cite{Rowe} for details on this subject.  
The method has also  been extended to Lie superalgebras \cite{Quesne}, 
and in recent years, to quantum groups and quantum supergroups
\cite{Lohe}.

A notable feature of the method of vector coherent states 
is that it is  particularly well adapted to standard techniques 
in physics,  thus is readily applicable to
addressing concrete physical problems.  
Apart from its physical applications, the method is also of great 
mathematical interests.  In particular, 
it is closely related \cite{Gilmore} to  
the Bott - Borel - Weil theorem.   
The theorem  is one of the hall marks in the representation 
theory of Lie groups.  It realizes   
the finite dimensional irreps of compact Lie groups 
in terms of  cohomology groups of homogeneous vector bundles.  
The Langlands and Kostant conjectures, proved by Schmid, 
generalize the theorem to noncompact semi - simple  Lie groups, 
yielding a geometrical realization of Harish - Chandra's 
discrete series of representations.   The Bott - Borel - Weil theory 
was further modified and significantly extended in the last 20 years, 
leading to the development of the theory of cohomological induction, 
which has now become a fundamental part of modern representation theory. 
The \BBW theory also plays important roles in various other fields, 
most notably, geometric quantization and Penrose transforms.  

The aim of the present note is to develop the vector coherent 
state construction of representations for affine Lie algebras.   
We will consider affine $sl(2)$ only in this note; 
more general cases will be treated else where. 
Affine Lie algebras \cite{Kac}  made their first appearance in physics 
in the context of dual resonance models,  now called strings ( See 
\cite{Goddard} for references on applications of affine Lie algebras 
to string theory.). 
Later studies showed that these algebras are also of paramount importance 
for conformal field theory and  for two dimensional 
integrable models in statistical mechanics(See \cite{Itzykson} 
for references.).  In fact the   
investigation of their applications in these  
three areas has been a major theme of mathematical physics since the 
early 1980s.  We should also mention that Wen proposed to 
utilize affine Lie algebras to describe the so - called 
edge states in fractional quantum Hall effect. 
We refer to Section 3 of Chapter 10 in Fradkin's book 
\cite{Fradkin}  for a concise discussion of Wen's proposal. 

We generalize the vector coherent state method to 
explicitly construct a large class of highest weight irreps, 
which include the integrable highest weight irreps as special 
cases. The main results are summarized in equations 
(\ref{simple}) and (\ref{general}), which we believe to be new. 
Our construction may be regarded as an infinitesimal version 
of the construction of Pressley and Segal's \cite{Pressley}, 
which yields integrable highest weight irreps of loop groups 
as holomorphic sections of vector bundles over 
`fundamental homogeneous spaces'.   
Our construction of $\g$ irreps also bears considerable
similarity with the  widely studied Wakimoto construction 
(See, e.g., \cite{Wakimoto}.), although it also 
differs from the latter in many essential ways. 
It will be interesting to find out what the exact 
relationship is between the two.   
Results presented here should extend to general affine 
Lie algebras and the Virasoro algebra, as well as their 
supersymmetric counter parts. We shall return to these 
problems in the future.         

\section{\normalsize AFFINE $sl(2)$ AND ITS REPRESENTATIONS} 
We will work on the complex number field $\C$. The affine Lie algebra 
$\g$ has a basis $\{\kappa; \ e[n],\ h[n], \ f[n], \ n\in\Z\}$ 
satisfying the following commutation relations
\ba 
{[} h[m], \ e[n] ] &=& 2\, e[m+n], \nonumber \\  
{[} h[m], \ f[n] ] &=&- 2\, f[m+n], \nonumber \\  
{[} h[m], \ h[n] ] &=& 2 m\,\kappa\, \delta_{m+n, 0},  \nonumber \\  
{[} e[m], \ f[n] ] &=& h[m+n] + m\,  \kappa\, \delta_{m+n, 0},  \nonumber \\  
\mbox{ the rest vanishes.} 
\na 
We will denote by $\U$ its universal enveloping algebra, 
and regard $\g$ as embedded in $\U$. Note that this affine Lie algebra 
contains a number of interesting subalgebras, which will be useful 
for studying the representation theory. In particular, we have 
\ban 
sl(2)&=& \langle e[0],\ h[0],\ f[0] \rangle,\\ 
{\frak b}_+&=&\langle \kappa;\ e[0],\ h[0];\, \ e[n],\ h[n], 
\ f[n], \ n>0\rangle, \\      
{\frak p}&=&\langle \kappa;;\ \ e[n],\ h[n], 
\ f[n], \ n\in\Z_+ \rangle,\\  
{\frak u}_+&=& \langle e[n],\ h[n], \ f[n], \ n>0 \rangle,\\ 
{\frak u}_-&=& \langle e[n],\ h[n], \ f[n], \ n<0 \rangle, 
\nan 
where ${\frak b}_+$ is a familiar Borel subalgebra, and  ${\frak p}$, 
containing  ${\frak b}_+$, will be called a parabolic subalgebra. 
Note also that $\g = {\frak u}_-  + {\frak p}$.

A standard technique in representation theory is to induce 
$\g$ modules, i.e., Verma modules, 
from given one dimensional ${\frak b}_+$ modules. 
More generally, one can also start with finite dimensional 
highest weight ${\frak p}$ modules to induce generalized Verma 
modules, and to obtain irreps as subquotients. 
This method is the analogue of the parabolic induction 
widely used in the representation theory of finite dimensional 
Lie algebras. It also provides the basis for developing the method 
of vector coherent states.

Let $V_0$ be a finite dimensional irreducible ${\frak p}$ module 
such that the positive modes $e[n],\ h[n], \ f[n], \ n>0$ act 
by zero.  This in particular implies the existence of a non - vanishing 
maximal vector $v_+\in V_0$ such that 
\ban 
e[0]\, v_+ &=&0,\\ 
h[0]\, v_+ &=& \lambda\,  v_+, \\ 
\kappa \, v_+ &=& c \, v_+, 
\nan    
where $\lambda\in{\Z}_+$, $c\in\C$.  Construct the vector space 
\ban 
W&=& \U\otimes_{U({\frak p})} V_0, 
\nan 
and define a left $\g$ action on it by the multiplication of $\U$, 
thus turning $W$ into a $\U$ module, which will be called a 
generalized Verma module.   It is clearly true that 
$W$ is $U({\frak u}_+)$ - locally  finite, i.e., 
for any $w\in W$ and $u\in {\frak u}_+$, 
there exists a finite positive integer 
$N$ such that $(u)^N\, w =0$.  Also, it is a direct consequence of 
the \BPW theorem that  $W=$ $U({\frak u}_-) \otimes_\C V_0$ 
as a vector space.  In general, $W$ itself is not 
 irreducible, but contains a unique maximal proper submodule $M$ 
such that the quotient 
\ban 
V&=& W/M 
\nan 
gives rise to an irreducible $\U$ module.  When $c-\lambda\in{\Z}_+$, 
$V$ is integrable,  and its structure is well understood \cite{Kac}.

From here on we will assume that $c\in{\Bbb R}$. This includes the 
integrable highest weight modules as special cases. 
Now we can give a relatively easy and explicit characterization of the 
maximal proper submodule $M$.  Recall that $\g$ admits the 
following anti - linear involution( denoted by $\dagger$ ) 
\ban      
     (\kappa)^\dagger &=& \kappa, \\  
     (e[n])^\dagger &=& f[-n], \\ 
     (f[n])^\dagger &=& e[-n],\\ 
     (h[n])^\dagger &=& h[-n], \ \ \ \ n\in\Z. 
\nan  
Using it, we can define a sesquilinear form $\langle \ | \ \rangle: $
$ W\times W \rightarrow \C$ by requiring
\ban 
&i).& \langle \ | \ \rangle|_{V_0\otimes V_0} \ \  
\mbox{defines a Hermitian inner product for $V_0$ regarded as}\\
& & \mbox{ a module of }  sl(2)=\langle e[0], \ h[0], \, f[0] \rangle, \ \
\mbox{and}\ \   \langle v_+ | v_+ \rangle =1;\\  
&ii). &  \langle u\, w_1 | w_2 \rangle 
= \langle w_1 | (u)^\dagger w_2 \rangle, 
\ \ \ \forall w_1, w_2\in W,  \ u\in\U.
\nan  
It is a straightforward exercise to show that $\langle \ | \ \rangle$ 
is well defined and unique. Furthermore, 
\ba 
M&= &Ker \langle \ | \ \rangle \\ 
&:=& \{ v\in W | \langle v |  W \rangle=0 \}. \nonumber    
\na

\section{\normalsize VECTOR COHERENT STATES} 
We now turn to the vector coherent state method, which, when 
restricted to integrable highest weight irreps, is the  
infinitesimal version of the construction of Pressley and 
Segal's \cite{Pressley}, which yields integrable highest weight 
irreps of loop groups as holomorphic sections of vector bundles over
`fundamental homogeneous spaces'.  Introduce 
the indeterminates $x_k,\, y_k,\, z_k, \ k=1, 2, ..., \infty$, 
and denote by $P(X, Y, Z)$ the linear span of their polynomials with 
coefficients in $\C$.  
Denote 
\ban 
E(X)&=&\sum_{k=1}^\infty e[k] x_k, \\       
H(Y)&=&\sum_{k=1}^\infty h[k] y_k, \\       
F(Z)&=&\sum_{k=1}^\infty f[k] z_k. 
\nan   
Introduce the formal power series,  
$\exp(E(X))$,  $\exp(H(Y)))$,  $\exp(F(Z))$,  
of the indeterminates with coefficients in $\U$, where, e.g.,  
\ban \exp(E(X)) &=&\sum_{n=0}^\infty (E(X))^n/n!.  \nan 
Let 
\ban 
g(X, Y, Z)&=&\exp(E(X))\, \exp(H(Y)))\, \exp(F(Z)).  
\nan  
Following the strategy of the vector coherent state method for 
ordinary Lie algebras \cite{Hecht, Rowe, Quesne, Gilmore}, 
we define the linear map $\xi(X, Y, Z): W\rightarrow 
P(X, Y, Z)\otimes_\C V_0$ by 
\ba 
\xi_w(X, Y, Z)&=&\sum_{i=1}^{dim V_0} \langle v_i | 
g(X, Y, Z) | w\rangle \otimes v_i, \ \ \ w\in W 
\na   
where $\{ v_i | i=1, 2, ..., dim V_0\}$ forms a basis of $V_0$ 
such that $\langle v_i | v_j \rangle=\delta_{i j}$, the existence 
of which is guaranteed by the finite dimensionality of $V_0$ 
as an $sl(2)$ module.   
We denote by ${\bf\Xi}(X, Y, Z)$ the image of the map, and call 
it the space of the vector coherent states.  Note that 
${\bf\Xi}(X, Y, Z)$ is indeed a subset of $P(X, Y, Z)\otimes_\C V_0$ 
because of the $U({\frak u}_+)$ - local finiteness of $W$. 

It is of crucial importance to observe that the coefficients of 
the formal power series $g(X, Y, Z)$ in the variables $x_k,\ y_k, \ z_k$, 
$k=1, 2, ..., \infty$, form a basis of $U({\frak u}_+)$.  Therefore, 
an element $w\in W$ belongs to the maximal proper submodule 
$M\subset W$ if and only if $\xi_w(X, Y, Z) =0$.  Also,  
${\bf\Xi}(X, Y, Z)$ admits a natural $\U$ action 
\ba 
u\cdot \xi_w(X, Y, Z)&=& \xi_{u w}(X, Y, Z), \ \ \ u\in\U, \ w\in W.
\na 
Now we have the following results:  
\ban  &i).&  Ker \xi(X, Y, Z) =M, \\ 
&ii).& {\bf\Xi}(X, Y, Z)\cong V \ \mbox{as irreducible}\  
\U \ \mbox{modules}.\\  \nan

One of the  advantages of the vector coherent state formulation 
of representations  is that the elements  of the affine 
Lie algebra can be realized explicit in terms of differential 
operators on $P(X, Y, Z)$ and endomorphisms of $V_0$.  
To obtain such a realization, 
we need some technical results, which we now discuss.  
Define the shifting operators $\sigma_n$, 
$n\in\Z$  by 
\ban 
\sigma_n E(X) &=& \sum_{k=1}^\infty e[n+k] x_k, \\ 
\sigma_n H(Y)&=& \sum_{k=1}^\infty h[n+k] y_k, \\
\sigma_n F(Z) &=& \sum_{k=1}^\infty f[n+k] z_k,  
\nan 
and also introduce the polynomials 
\ban 
{\bf Z}_N(Y)&=& { {1}\over{N!} } { {d^N}\over{d t^N} } 
    \exp(\sum_{k=1}^\infty t^k y_k)\left\arrowvert_{t=0}\right., 
\ \ \ N\in\Z_+. 
\nan   
Some lengthy, but relatively straightforward, calculations yield 
the following relations,  
\ba {[} \exp(E(X)),\ h[k]\, ] &=& - 2 \exp(E(X)) \sigma_k E(X),\nonumber \\
{[} \exp(H(Y))),\ h[k]\, ] &=&  
                   -2\kappa k \theta(-k)  y_{-k} \exp(H(Y))),\nonumber \\  
{[} \exp(F(Z)),\ h[k]\, ] &=&  2 \exp(F(Z)) \sigma_k F(Z),\nonumber \\
f[k] \exp(H(Y))) &=& \exp(H(Y))) \sum_{N=0}^{\infty} 
                 {\bf Z}_N(2 Y) f[k+N], \nonumber\\ 
{[} \exp(E(X)),\ f[k]\, ] &=&  \exp(E(X))\left\{  \sigma_k H(X) 
            - k \kappa x_{-k} \theta(- k) + 
            \sum_{p=1}^\infty  x_p \sigma_{p+k}E(X) \right\} , \nonumber\\ 
{[} \exp(F(Z)),\ e[k]\, ]&=& \exp(F(Z))\left\{ \sum_{p=1}^\infty
 z_p \sigma_{p+k} F(Z) - \kappa k z_{-k} \theta(- k) - \sigma_k H(Z)
\right\},\nonumber \\  
\exp(H(Y)) e[k]&=& \sum_{N=0}^{\infty} {\bf Z}_N(2 Y) e[k+N] \exp(H(Y)), 
\ \ \ k\in\Z,  
\na 
where $2 Y= (2 y_1, 2 y_2, ...)$,  and  
 $\theta(k) =\left\{\begin{array}{l l}
                        1,  & k> 0,\\
                        0, & k\le0.  
                      \end{array}\right. $ 

With the help of these relations, we can readily obtain the 
realization of the elements of the affine Lie algebra $\g$ in 
explicit form.  We will denote the realization by 
$\xi$.  For the simple and Cartan generators, we have 
\ba 
\xi(\kappa)&=& c, \nonumber\\ 
\xi(f[1])&=&{{\partial}\over{\partial z_1}},\nonumber \\
\xi(h[0])&=& \pi_0(h) + 2 \sum_{p=1}^\infty 
    \left[ z_p {{\partial}\over{\partial z_p}} - 
    x_p {{\partial}\over{\partial x_p}}\right], \nonumber \\   
\xi(e[0])&=& \pi_0(e) + \sum_{N=1}^\infty {\bf Z}_N( 2 Y)
    {{\partial}\over{\partial x_N}} - \sum_{p=1}^\infty \left[ 
    z_p {{\partial}\over{\partial y_p}} + z_p\sum_{q=1}^\infty 
    z_q{{\partial}\over{\partial z_{p+q}}}\right],\nonumber \\
\xi(f[0])&=& \pi_0(f) - \sum_{N=1}^\infty {\bf Z}_N( 2 Y)
    {{\partial}\over{\partial z_N}}  
    + \sum_{p=1}^\infty \left[ 
      x_p {{\partial}\over{\partial y_p}} +
    x_p \sum_{q=1}^{\infty} x_q {{\partial}\over{\partial x_{p+q}}}\right], 
    \nonumber\\ 
\xi(e[-1])&=&2\, y_1 \pi_0(e) + z_1\left\{ c - \pi_0(h) + 
      2 \sum_{p=1}^\infty x_p {{\partial}\over{\partial x_p}} 
      \right\}\nonumber\\ 
    & +& \sum_{N=1}^\infty {\bf Z}_{N+1}( 2 Y) 
     {{\partial}\over{\partial x_N}}  
   -\sum_{p=1}^\infty \left\{  z_{p+1}{{\partial}\over{\partial y_p}} 
   + z_p  \sum_{q=1}^\infty z_q {{\partial}\over{\partial z_{p+q-1} } }  
      \right\}   ,  
\label{simple} 
\na  
where $\pi_0(e), \, \pi_0(h),\, \pi_0(f)$ $\in End(V_0)$ are respectively 
defined by the actions of $e[0],\, h[0]$ and $f[0]$ on $V_0$.  

In principle, these elements are sufficient to generate the entire 
$\g$ algebra. However, it is not much more difficult to work out the 
realization for other generators, which we spell out below:
\ba 
\xi(f[k])&=&{{\partial}\over{\partial z_k}},\nonumber \\  
\xi(h[k])&=& {{\partial}\over{\partial y_k} } + 2 \sum_{p=1}^\infty
      z_p {{\partial}\over{\partial z_{k+p} }},\nonumber \\ 
\xi(e[k])&=& {{\partial}\over{\partial x_k}} 
      + \sum_{N=1}^\infty {\bf Z}_{N}( 2 Y) 
      {{\partial}\over{\partial x_{k+N}}}  
      - \sum_{p=1}^\infty \left[
    z_p {{\partial}\over{\partial y_{k+p}}} + z_p\sum_{q=1}^\infty
    z_q{{\partial}\over{\partial z_{k+p+q}}}\right],\nonumber \\   
\xi(f[-k])&=&{\bf Z}_k (-2 Y) \xi(f[0]) 
           +\sum_{t=1}^k {\bf Z}_{k-t}(-2 Y) {\cal D}_t, \nonumber \\
\xi(h[-k])&=& 2 c k y_k - 2 x_k \pi_0(e) 
      + 2\sum_{p=1}^\infty\left[z_{k+p}{{\partial}\over{\partial z_p}}
        - x_{k+p} {{\partial}\over{\partial x_p}}\right] 
       +2 x_k \sum_{p=1}^k z_p \xi(f[p-k]), \nonumber \\ 
\xi(e[-k])&=& c k z_k + {\bf Z}_k(2 Y) \pi_0(e) 
     + \sum_{N=1}^\infty {\bf Z}_{N+k}
     (2 Y) {{\partial}\over{\partial x_N}} \nonumber \\
     &-&\xi( \sigma_{-k}H(Z)) + \sum_{p=1}^\infty z_p\xi(\sigma_{p-k} F(Z)),  
\hspace{3.5cm} k>0, 
\label{general}     
\na 
where 
\ban 
{\cal D}_k&=& x_k ( c k + \pi_0(h) ) \\ 
          &-& \sum_{N=1}^\infty {\bf Z}_{N+k}( 2 Y) 
               {{\partial}\over{\partial z_N}} 
           +\sum_{p=1}^\infty x_{p+k}\left( {{\partial}\over{\partial y_p}} 
           +\sum_{q=1}^\infty x_q {{\partial}\over{\partial x_{p+q}}}\right)\\
          &-&\sum_{p=1}^\infty \sum_{q=1}^k x_p x_q 
          \left\{ \theta(p+q-k) {{\partial}\over{\partial x_{p+q-k}}}
           +\delta_{p+q, k}\pi_0(e)\right\}, \ \ \ \ \ \ k>0.
\nan     
We can also introduce the grading operator $d$. It necessarily acts 
on $V_0$ by a constant $d_0\in{\Bbb R}$, which can in fact be 
chosen arbitrarily.  Now    
\ban 
\xi(d)&=&d_0 - \sum_{p=1}^\infty p \left( x_p {{\partial}\over{\partial  x_p}}
          + y_p {{\partial}\over{\partial y_p}} 
          + z_p {{\partial}\over{\partial z_p}} \right). \\    
\nan

An important feature of the vector coherent state realization is that 
repeated applications of  $\xi(e[-k])$, $\xi(h[-k])$, $\xi(f[-k])$, 
$k>0$, to $1\otimes V_0$ will automatically generate the irreducible 
$\U$ module ${\bf\Xi}(X, Y, Z)$. 
Also observe that $\xi(\g)$ acts naturally on ${\bar W}(X, Y, Z)=$ 
$P(X, Y, Z)\otimes V_0$, and ${\bf\Xi}(X, Y, Z)$ is contained 
in ${\bar W}(X, Y, Z)$ as the unique minimal submodule.  
It appears possible to extend the K - matrix technique\cite{Rowe}
to the present context to explicitly identify this submodule. 
A different approach is to find, by cohomological means, 
a system of linear differential operators 
with constant coefficients acting on ${\bar W}(X, Y, Z)$, 
the kernel of which will reproduce ${\bf\Xi}(X, Y, Z)$. 
We will report our results on 
these problems in a separate publication, where the coherent 
state method is also developed for general affine Lie algebras.

\end{document}